# DEVELOPMENT OF THREE DIMENSIONAL CONSTITUTIVE THEORIES BASED ON LOWER DIMENSIONAL EXPERIMENTAL DATA

Satish Karra, K. R. Rajagopal, College Station

*Abstract.* Most three dimensional constitutive relations that have been developed to describe the behavior of bodies are correlated against one dimensional and two dimensional experiments. What is usually lost sight of is the fact that infinity of such three dimensional models may be able to explain these experiments that are lower dimensional. Recently, the notion of maximization of the rate of entropy production has been used to obtain constitutive relations based on the choice of the stored energy and rate of entropy production, etc. In this paper we show different choices for the manner in which the body stores energy and dissipates energy and satisfies the requirement of maximization of the rate of entropy production that leads to many three dimensional models. All of these models, in one dimension, reduce to the model proposed by Burgers to describe the viscoelastic behavior of bodies.

## 1. Introduction

An observation of a phenomenon or a set of phenomena leads one to conjecture as to its cause and forms the basis for the crude first step in the development of a model. An experiment is then deliberately and carefully fashioned to test and refine the conjecture. Unfortunately, this procedure is rendered most daunting as one is not usually accorded the luxury of being able to perform sufficiently general three dimensional experiments while the models that one would like to develop are fully three dimensional models. Most of the general three dimensional constitutive models that are being used in continuum mechanics have been developed on the basis of information gleaned from one or two dimensional special experiments. It does not take much mathematical acumen to recognize the dangers fraught in the process of such generalizations as infinity of three dimensional models could be capable of explaining the lower dimensional experimental data. Of course, one does not corroborate a three dimensional model by merely comparing against data from a single one dimensional experiment. One tests the model against several different experiments, but these experiments tend to be simple experiments in view of the extraordinary difficulties in developing an experimental program that can truly test the full three dimensionality of the model, especially when the response that is being described is complex. In order to obtain a meaningful three dimensional model on the basis of experimental data in lower dimensions, one needs to be guided by enormous physical insight and intuition. This is easier said than done and in fact most models that are currently in use are based on flimsy and tenuous rationale.

One might be tempted to think that the dictates of physics would greatly aid in the development of models from experimental data. For instance, the second law of thermodynamics could play a stringent role in the class of admissible models. Similarly, invariance requirements such as Galilean invariance could also provide restrictions on the class of allowable models. Unfortunately, the sieve provided by all such restrictions is far too coarse as it permits several models to go through that exhibit undesirable properties.

While modeling, one might start directly by assuming a constitutive relation between the stress and other relevant quantities. This relation could be an explicit expression (function) for the stress in terms of kinematical variables as in the case of Hooke's law or the Navier-Stokes model, or it could be an implicit relation as in the case of many rate type non-Newtonian fluid models. Assuming such constitutive relations implies six scalar constitutive relations (in the case of the stress being symmetric). One could also assume forms for the manner in which energy is stored and entropy is produced by the body and determine the constitutive relation for the stress by appealing to a general thermodynamic framework that has been put in place (see the review articles by Rajagopal and Srinivasa [1], [2] for details of the framework). The framework casts the second law as an equation that defines the rate of entropy production (see Green and Naghdi [3], Rajagopal and Srinivasa [4]) and appeals to the maximization of the rate of entropy production (while Ziegler [5] had appealed to such a requirement, the context within which he made such an appeal is different from that required by Rajagopal and Srinivasa [1], [2]). The general thermodynamic framework has been used to describe a plethora of disparate material response: viscoelasticity, inelasticity, twinning, phase transition in solids, behavior of single crystals super alloys, mixtures, inhomogeneous fluid, etc. While the method seems exceedingly powerful, there are some interesting nuances concerning its application that the modeler should be aware of, and in this paper by constructing explicit examples we illustrate these delicate issues. It is important to recognize that one can obtain the same constitutive relations for the stress by choosing different forms for the stored energy functions and the rate of entropy production (see Rao and Rajagopal [6] who develop the non-linear three dimensional Maxwell model by choosing two different sets of stored energy and rate of dissipation). In fact, it is possible that several sets of stored energy and rate of dissipation function can lead to the same model. We illustrate this by considering four different sets of stored energy and rate of dissipation to obtain the model developed by Burgers [7], and these four choices are different from a previous choice made by Murali Krishnan and Rajagopal [8]. It is interesting to note that by making the choice of two scalar functions, we can arrive at a constitutive relation for the stress, a tensor with six scalar components. Many of the one-dimensional models that have been developed to describe the response of viscoelastic materials was by appealing to an analogy to mechanical systems of springs (means for storing energy), and dashpots (means for dissipating energy/ producing entropy), though in his seminal paper on viscoelasticity Maxwell [9] did not appeal to such an analogy. Within the context of these mechanical systems, it becomes clear how one can get the same form for the stress by choosing different stored energy and rate of entropy production functions as one can choose different networks of springs and dashpots to effect the same response.

In 1934 Burgers [7] developed the following one-dimensional model by appealing to a mechanical analog:

$$(1.1) \qquad \sigma + p_1 \dot{\sigma} + p_2 \ddot{\sigma} = q_1 \dot{\epsilon} + q_2 \ddot{\epsilon},$$

where $p_1$ and $p_2$ are relaxation times, $q_1$ and $q_2$ are viscosities, and $\sigma$ and $\epsilon$ denote the stress and the linearized strain respectively. A three dimensional generalization of that was provided by Murali Krishnan and Rajagopal [8], within the context of a thermodynamic basis that requires that among an admissible class of constitutive relations that which is selected is the one that maximizes the rate of entropy production. The second law merely requires that the entropy production be non-negative and one would expect the requirement of maximization of the rate of entropy to cull the class of rate of entropy production functions. As we shall restrict our analysis to a purely mechanical context, instead of making a choice for the rate of entropy production we shall make a choice for the rate of dissipation (the rate at which working is converted to heat) which is the only way in which entropy is produced within the context of interest.

We shall assume that the class of bodies we are interested in modeling are viscoelastic fluids that are capable of instantaneous elastic response. A body that exists in a configuration $\kappa_t$ under the action of external stimuli, on the removal of the external stimuli could attain a configuration $\kappa_{p(t)}$, which is referred to as a natural configuration corresponding to the configuration $\kappa_t$. However, more than one natural configuration could be associated with the configuration $\kappa_t$ based on how the external stimuli is removed, whether instantaneously, slowly, etc. The natural configuration that is accessed depends on the process class allowed. In this study, we shall assume the natural configuration that is achieved is that due to an instantaneous unloading to which the body responds in an instantaneous elastic manner. A detailed discussion of the role of natural configurations can be found in Rajagopal [10] and the review article by Rajagopal and Srinivasa [1]. Even within the context of an instantaneous elastic unloading, it might be possible that the body could go to different natural configurations $\kappa_{p_i(t)}$, $i = 1, ...., n$.

When one provides a spring-dashpot mechanical analogy for a viscoelastic material one obtains a constitutive equation that holds at a point, i.e., the point is capable of storing energies like the various springs and dissipate energy as the dashpots, but it also has to take into account the arrangement of the springs and the dashpots. The central idea of Mixture Theory is that the various constituents of the mixture co-exist. That is, in a homogenized sense at a point, the model has to reflect the combined storage of energies in the springs and the dissipation of the dashpots based on the way in which they are arranged. Of course, a point is a mathematical creation that does not exist, and what is being modeled is a sufficiently small chunk in the body. This chunk can store and dissipate energy in different ways. The point of importance is various arrangements of springs and dashpots can lead to the same net storage of energy of the springs and the dissipation by the dashpots. Put differently, the chunk can respond in an identical manner for different ways in which the springs and dashpots are put together. This is essentially the crux of the paper. We have five different three dimensional models, four that are developed in this paper and one that was developed by Murali Krishnan and Rajagopal [8] and all five three dimensional models could claim equal status as generalizations of the one dimensional model developed by Burgers. Recently, Malek and Rajagopal [11] used the thermodynamic framework that we have discussed to obtain a model for two viscous liquids. In this paper, we have a more complicated mixture in that we have two different elastic solids coexisting with two different dissipative fluids. We do not allow for relative motion between the constituents, we assume they coexist and move together.

The organization of the paper is as follows. In section 2, we introduce the kinematics that is necessary to the study and the basic balance laws for mass, linear and angular momentum. We also introduce the second law of thermodynamics. This introduction is followed by a discussion of four different models which all reduce to Burgers' one dimensional model in sections 3–6. We make a few final remarks in the last section.

## 2. Preliminaries

Let $\kappa_R(\mathcal{B})$ and $\kappa_t(\mathcal{B})$ denote, respectively the reference configuration of the body, and the configuration of the body $\mathcal{B}$ at time $t$. Let $\mathbf{X}$ denote a typical point belonging to $\kappa_R(\mathcal{B})$ and $\mathbf{x}$ the same material point at time $t$, belonging to $\kappa_t(\mathcal{B})$. Let $\chi_{\kappa_R}$ denote a sufficiently smooth mapping that assigns to each $\mathbf{X} \in \kappa_R(\mathcal{B})$, a point $\mathbf{x} \in \kappa_t(\mathcal{B})$, i.e.,

$$\mathbf{x} := \chi_{\kappa_R}(\mathbf{X}, t). \tag{2.1}$$

The velocity $\mathbf{v}$, the velocity gradient $\mathbf{L}$ and the deformation gradient $\mathbf{F}_{\kappa_R}$ are defined through

$$\mathbf{v} := \frac{\partial \chi_{\kappa_R}}{\partial t}, \quad \mathbf{L} := \frac{\partial \mathbf{v}}{\partial \mathbf{x}}, \quad \mathbf{F}_{\kappa_R} := \frac{\partial \chi_{\kappa_R}}{\partial \mathbf{X}}. \tag{2.2}$$

It immediately follows that

$$\mathbf{L} = \dot{\mathbf{F}}_{\kappa_R} \mathbf{F}_{\kappa_R}^{-1}. \tag{2.3}$$

We denote the symmetric part of the velocity gradient by $\mathbf{D}$, i.e.,

$$\mathbf{D} := \frac{1}{2}\left(\mathbf{L} + \mathbf{L}^T\right). \tag{2.4}$$

The left and right Cauchy-Green stretch tensors $\mathbf{B}_{\kappa_R}$ and $\mathbf{C}_{\kappa_R}$ are defined through

$$\mathbf{B}_{\kappa_R} := \mathbf{F}_{\kappa_R} \mathbf{F}_{\kappa_R}^T, \quad \mathbf{C}_{\kappa_R} := \mathbf{F}_{\kappa_R}^T \mathbf{F}_{\kappa_R}. \tag{2.5}$$

Let $\kappa_{p(t)}$ denote the preferred natural configuration associated with the configuration $\kappa_t$. We define $\mathbf{F}_{\kappa_{p(t)}}$ as the mapping from the tangent space at a material point in $\kappa_{p(t)}$ to the tangent space at the same material point at $\kappa_t$ (see Fig. 1). We then define

$$\mathbf{B}_{\kappa_{p(t)}} := \mathbf{F}_{\kappa_{p(t)}} \mathbf{F}_{\kappa_{p(t)}}^T, \quad \mathbf{C}_{\kappa_{p(t)}} := \mathbf{F}_{\kappa_{p(t)}}^T \mathbf{F}_{\kappa_{p(t)}}. \tag{2.6}$$

The mapping $\mathbf{G}$ is defined through (see Fig. 1)

$$\mathbf{G} := \mathbf{F}_{\kappa_R \to \kappa_{p(t)}} := \mathbf{F}_{\kappa_{p(t)}}^{-1} \mathbf{F}_{\kappa_R}. \tag{2.7}$$

We can then define the tensor $\mathbf{C}_{\kappa_R \to \kappa_{p(t)}}$ in a manner analogous to $\mathbf{C}_{\kappa_R}$ through

$$\mathbf{C}_{\kappa_R \to \kappa_{p(t)}} := \mathbf{G}^T \mathbf{G}, \tag{2.8}$$

and it follows that

$$\mathbf{B}_{\kappa_{p(t)}} = \mathbf{F}_{\kappa_R} \mathbf{C}_{\kappa_R \to \kappa_{p(t)}}^{-1} \mathbf{F}_{\kappa_R}^T. \tag{2.9}$$

We shall also record balance of mass (assuming incompressibity), balance of linear and angular momentum (in the absence of body couples):

$$div(\mathbf{v}) = 0, \quad \rho\dot{\mathbf{v}} = div(\mathbf{T}^T) + \rho\mathbf{b}, \quad \mathbf{T} = \mathbf{T}^T, \tag{2.10}$$

where $\rho$ is the density, $\mathbf{v}$ is the velocity, $\mathbf{T}$ is the Cauchy stress tensor, $\mathbf{b}$ is the specific body force, $div(.)$ is the divergence operator with respect to the current configuration and $(.)^T$ denotes transpose. In addition, the local form of balance of energy is

$$\rho\dot{\epsilon} = \mathbf{T}.\mathbf{L} - div(\mathbf{q}) + \rho r, \tag{2.11}$$

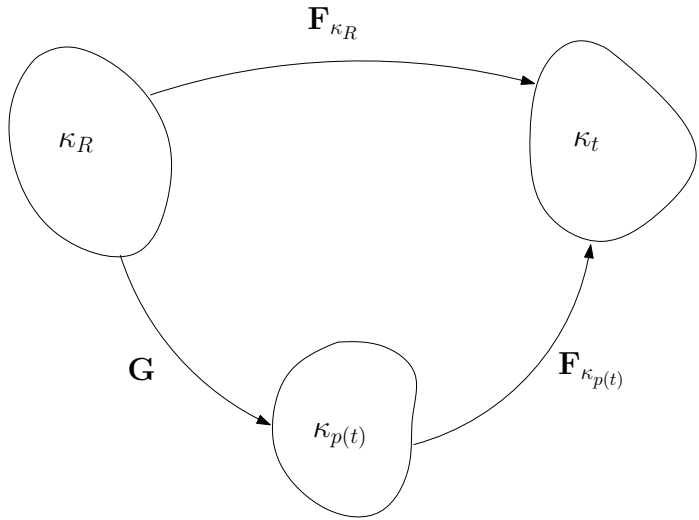

FIGURE 1. Schematic of the natural configuration $\kappa_{p(t)}$ corresponding to the current configuration $\kappa_t$ and the relevant mappings from the tangent spaces of the same material point in $\kappa_R$, $\kappa_t$ and $\kappa_{p(t)}$.

where $\epsilon$ denotes the specific internal energy, $\mathbf{q}$ denotes the heat flux vector and $r$ denotes the specific radiant heating. We shall invoke the second law of thermodynamics in the form of the *reduced energy dissipation equation*, for isothermal processes:

(2.12) $$\mathbf{T}.\mathbf{D} - \rho\dot{\psi}|_{isothermal} := \xi \geq 0,$$

where $\psi$ is the specific Helmholtz potential, $\xi$ denotes the rate of dissipation (specifically rate of entropy production).

When one works with implicit constitutive models of the form

(2.13) $$\mathbf{f}(\mathbf{T}, \mathbf{D}) = 0,$$

where $\mathbf{T}$ is the Cauchy stress, or more general models of the form

(2.14) $$\mathbf{f}\left(\mathbf{T}, \overset{\nabla}{\mathbf{T}}, \ldots, \overset{(n)}{\overset{\nabla}{\mathbf{T}}}, \mathbf{D}, \overset{\nabla}{\mathbf{D}}, \ldots, \overset{(n)}{\overset{\nabla}{\mathbf{D}}}\right) = 0,$$

where the superscript $\overset{(n)}{\nabla}$ stands for the $n$ Oldroyd derivatives [12], and where $\mathbf{T}$ and $\mathbf{D}$ seem to have the same primacy in that the maximization could be with respect to $\mathbf{T}$ or $\mathbf{D}$. However, the superficial assumption that $\mathbf{T}$ and $\mathbf{D}$ have the same primacy is incorrect as $\mathbf{T}$ (or the applied traction which leads to the stresses) causes the deformation (the appropriate kinematic tensor). In order to get sensible constitutive equations one ought to keep $\mathbf{D}$ fixed and vary $\mathbf{T}$ to find how a body responds to the stress that is a consequence of the applied traction. This comes up naturally in the development of implicit constitutive theories (see Rajagopal and Srinivasa [13], Rajagopal and Srinivasa [14]). More recently, Rajagopal [15] has discussed at length the implicit nature of constitutive relations. When one thinks explicity along classical terms of the stress being given explicity in terms of the kinematical variables, it is natural to hold $\mathbf{T}$ fixed and maximize with respect to the kinematical variable, in our case $\mathbf{D}$. This is what is followed in this work.

# 3. Model 1

3.1. **Preliminaries.** Let $\kappa_R$ denote the undeformed reference configuration of the body. We shall assume that the body has associated with it two natural configurations i.e., configurations to which it can be instantaneously elastically unloaded and corresponds to two mechanisms for storing energy (within one dimensional mechanical analog – two springs). Interestingly, one can get from the reference configuration to the two evolving natural configurations denoted by $\kappa_{p_i(t)}$, $i = 1, 2$ (see Fig. 2), via two dissipative responses. Let $\mathbf{F}_i$, $i = 1, 2, 3$, denote the gradients of the motion[1] from $\kappa_R$ to $\kappa_{p_1(t)}$, $\kappa_{p_1(t)}$ to $\kappa_{p_2(t)}$, and $\kappa_{p_2(t)}$ to $\kappa_t$ respectively. Also, we shall define the left Cauchy-Green stretch tensors,

$$\mathbf{B}_i := \mathbf{F}_i \mathbf{F}_i^T, \quad i = 1, 2, 3, \tag{3.1}$$

and the velocity gradients with their corresponding symmetric parts,

$$\mathbf{L}_i := \dot{\mathbf{F}}_i \mathbf{F}_i^{-1}, \quad \mathbf{D}_i := \frac{1}{2}\left(\mathbf{L}_i + \mathbf{L}_i^T\right), \quad i = 1, 2, 3. \tag{3.2}$$

Also, we note that [2]

$$\mathbf{F} = \mathbf{F}_3 \mathbf{F}_2 \mathbf{F}_1. \tag{3.3}$$

Let us denote the gradient of the motion from $\kappa_{p_1(t)}$ to $\kappa_t$ by $\mathbf{F}_p$; then,

$$\mathbf{F}_p = \mathbf{F}_3 \mathbf{F}_2, \tag{3.4}$$

and

$$\mathbf{F} = \mathbf{F}_p \mathbf{F}_1. \tag{3.5}$$

The left Cauchy-Green stretch tensor, the velocity gradient with its symmetric part, corresponding to $\mathbf{F}_p$ are

$$\mathbf{B}_p := \mathbf{F}_p \mathbf{F}_p^T, \quad \mathbf{L}_p := \dot{\mathbf{F}}_p \mathbf{F}_p^{-1}, \quad \mathbf{D}_p := \frac{1}{2}\left(\mathbf{L}_p + \mathbf{L}_p^T\right), \tag{3.6}$$

respectively.

Now, taking the time derivative of eqn. (3.5) we get:

$$\begin{aligned} \dot{\mathbf{F}} &= \dot{\mathbf{F}}_p \mathbf{F}_1 + \mathbf{F}_p \dot{\mathbf{F}}_1 \\ \Rightarrow \quad \mathbf{L}\mathbf{F} &= \mathbf{L}_p \mathbf{F}_p \mathbf{F}_1 + \mathbf{F}_p \mathbf{L}_1 \mathbf{F}_1 \\ \Rightarrow \quad \mathbf{L} &= \mathbf{L}_p + \mathbf{F}_p \mathbf{L}_1 \mathbf{F}_p^{-1}. \end{aligned} \tag{3.7}$$

Similarly, taking the time derivative of eqn. (3.4), we arrive at

$$\mathbf{L}_p = \mathbf{L}_3 + \mathbf{F}_3 \mathbf{L}_2 \mathbf{F}_3^{-1}. \tag{3.8}$$

Now,

$$\begin{aligned} \dot{\mathbf{B}}_p &= \mathbf{F}_p \dot{\mathbf{F}}_p^T + \dot{\mathbf{F}}_p \mathbf{F}_p^T \\ &= \mathbf{F}_p \mathbf{F}_p^T \mathbf{L}_p^T + \mathbf{L}_p \mathbf{F}_p \mathbf{F}_p^T \\ &= \mathbf{B}_p \mathbf{L}_p^T + \mathbf{L}_p \mathbf{B}_p. \end{aligned} \tag{3.9}$$

---

[1]In general, these are appropriate mappings of tangent spaces containing the same material point in different configurations.

[2]Henceforth, we shall denote $\mathbf{F}_{\kappa_R}$ by $\mathbf{F}$.

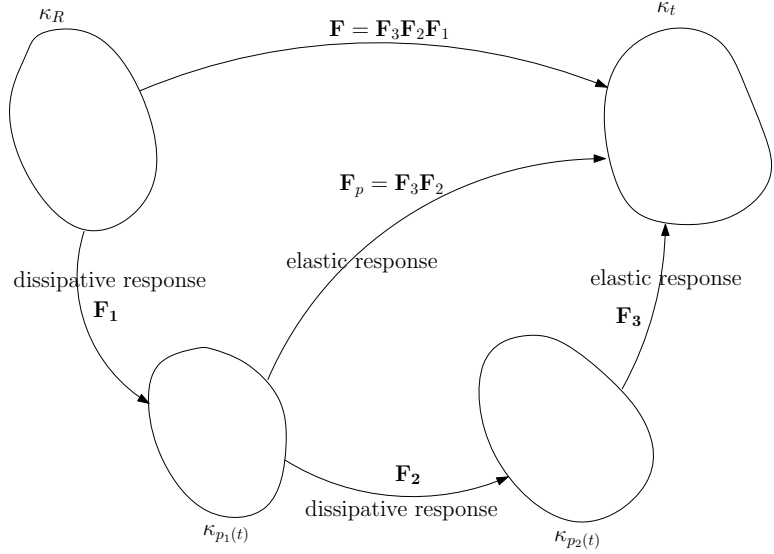

FIGURE 2. Schematic to illustrate the natural configurations for model 1. $\kappa_R$ is the reference configuration, $\kappa_t$ denotes the current configuration, and $\kappa_{p_1(t)}, \kappa_{p_2(t)}$ denote the two evolving natural configurations. The body dissipates energy like a viscous fluid as it moves from, $\kappa_R$ to $\kappa_{p_1(t)}$, and $\kappa_{p_1(t)}$ to $\kappa_{p_2(t)}$. Also, as shown, the body stores energy during its motion from, $\kappa_{p_2(t)}$ to $\kappa_t$, and $\kappa_{p_1(t)}$ to $\kappa_t$.

Post-multiplying eqn. (3.7) by $\mathbf{B}_p$, pre-multiplying the transpose of eqn. (3.7) by $\mathbf{B}_p$, and adding, we obtain

$$\stackrel{\nabla}{\mathbf{B}}_p = -2\mathbf{F}_p \mathbf{D}_1 \mathbf{F}_p^T, \tag{3.10}$$

where $\stackrel{\nabla}{\mathbf{B}}_p := \dot{\mathbf{B}}_p - \mathbf{B}_p \mathbf{L}^T - \mathbf{L}\mathbf{B}_p$ is the Oldroyd derivative of $\mathbf{B}_p$. In a similar fashion, using eqn. (3.8) and the relation $\dot{\mathbf{B}}_3 = \mathbf{B}_3 \mathbf{L}_3^T + \mathbf{L}_3 \mathbf{B}_3$, we get

$$\stackrel{\nabla_p}{\mathbf{B}}_3 = -2\mathbf{F}_3 \mathbf{D}_2 \mathbf{F}_3^T, \tag{3.11}$$

where $\stackrel{\nabla_p}{\mathbf{B}}_3 := \dot{\mathbf{B}}_3 - \mathbf{B}_3 \mathbf{L}_p^T - \mathbf{L}_p \mathbf{B}_3$. This is same as the Oldroyd derivative of $\mathbf{B}_3$, when the natural configuration $\kappa_{p_1(t)}$ is made the reference configuration.

We also note from eqn. (3.3) that

$$\begin{aligned}\dot{\mathbf{F}} &= \dot{\mathbf{F}}_3 \mathbf{F}_2 \mathbf{F}_1 + \mathbf{F}_3 \dot{\mathbf{F}}_2 \mathbf{F}_1 + \mathbf{F}_3 \mathbf{F}_2 \dot{\mathbf{F}}_1 \\ \Rightarrow \mathbf{L} &= \mathbf{L}_3 + \mathbf{F}_3 \mathbf{L}_2 \mathbf{F}_3^{-1} + \mathbf{F}_3 \mathbf{F}_2 \mathbf{L}_1 \mathbf{F}_2^{-1} \mathbf{F}_3^{-1}.\end{aligned} \tag{3.12}$$

and hence

$$\begin{aligned}\mathbf{I}.\dot{\mathbf{B}}_3 &= \mathbf{I}.\left(\mathbf{L}_3 \mathbf{B}_3 + \mathbf{B}_3 \mathbf{L}_3^T\right) \\ &= \mathbf{I}.\left(\mathbf{L}\mathbf{B}_3 - \mathbf{F}_3 \mathbf{L}_2 \mathbf{F}_3^T - \mathbf{F}_3 \mathbf{F}_2 \mathbf{L}_1 \mathbf{F}_2^{-1}\mathbf{F}_3^T + \mathbf{B}_3 \mathbf{L}^T - \mathbf{F}_3 \mathbf{L}_2^T \mathbf{F}_3^T - \mathbf{F}_3 \mathbf{F}_2^{-T} \mathbf{L}_1^T \mathbf{F}_2^T \mathbf{F}_3^T\right) \\ &= 2\mathbf{B}_3.\mathbf{D} - 2\mathbf{C}_3.\mathbf{D}_2 - \mathbf{C}_3.\left(\mathbf{F}_2 \mathbf{L}_1 \mathbf{F}_2^{-1} + \mathbf{F}_2^{-T} \mathbf{L}_1^T \mathbf{F}_2^T\right).\end{aligned} \tag{3.13}$$

The relations derived, in this sub-section, are sufficient for the purpose of analysing this model. In the following sub-section, we shall constitutively specify the forms for storage and rate of dissipation functions, and then we shall maximize the rate of dissipation subject to appropriate constraints (incompressibility and the energy dissipation equation), to determine the constitutive relation.

3.2. **Constitutive assumptions.** Let us assume the specific stored energy $\psi$ and the rate of dissipation $\xi$ of the form[3]

$$\psi \equiv \psi(\mathbf{B}_3, \mathbf{B}_p), \quad \xi \equiv \xi(\mathbf{D}_1, \mathbf{D}_2). \tag{3.14}$$

In particular, assuming that the instantaneous elastic responses from $\kappa_{p_1(t)}$ and $\kappa_{p_1(t)}$ are isotropic, and in virtue of incompressibility of the body, we choose

$$\psi(\mathbf{B}_3, \mathbf{B}_p) = \frac{\mu_3}{2\rho}(\mathbf{I}.\mathbf{B}_3 - 3) + \frac{\mu_p}{2\rho}(\mathbf{I}.\mathbf{B}_p - 3), \tag{3.15}$$

and

$$\xi(\mathbf{D}_1, \mathbf{D}_2) = \eta_1' \mathbf{D}_1.\mathbf{D}_1 + \eta_2' \mathbf{D}_2.\mathbf{D}_2. \tag{3.16}$$

The above assumption means that the body possesses instantaneous elastic response from the two evolving natural configurations $(\kappa_{p_1(t)}, \kappa_{p_2(t)})$ to the current configuration $\kappa_t$ (Fig. 2); the body stores energy like a neo-Hookean solid during its motion, from $\kappa_{p_1(t)}$ to $\kappa_t$, and from $\kappa_{p_2(t)}$ to $\kappa_t$. In addition, the response is linear viscous fluid-like, as the body moves from $\kappa_R$ to $\kappa_{p_1(t)}$, and from one natural configuration $(\kappa_{p_1(t)})$ to the other $(\kappa_{p_2(t)})$.

Also, since we have assumed that the material's instantaneous elastic response is isotropic, we shall choose the configurations $\kappa_{p_1(t)}, \kappa_{p_2(t)}$ such that

$$\mathbf{F}_3 = \mathbf{V}_3, \quad \mathbf{F}_p = \mathbf{V}_p, \tag{3.17}$$

where $\mathbf{V}_3, \mathbf{V}_p$ are the right stretch tensors in the polar decomposition i.e., the natural configurations are appropriately rotated.

Finally, using eqns. (3.13) and (3.17), we get

$$\mathbf{I}.\dot{\mathbf{B}}_3 = 2\mathbf{B}_3.\left[\mathbf{D} - \mathbf{D}_2 - \frac{1}{2}\left(\mathbf{F}_2 \mathbf{L}_1 \mathbf{F}_2^{-1} + \mathbf{F}_2^{-T} \mathbf{L}_1^T \mathbf{F}_2^T\right)\right], \tag{3.18}$$

and similarly

$$\mathbf{I}.\dot{\mathbf{B}}_p = 2\mathbf{B}_p.(\mathbf{D} - \mathbf{D}_1). \tag{3.19}$$

Substituting eqns. (3.15), (3.16) into (2.12) and using the relations in eqns. (3.18), (3.19),

$$\mathbf{T}.\mathbf{D} - \mu_3 \mathbf{B}_3.\left[\mathbf{D} - \mathbf{D}_2 - \frac{1}{2}\left(\mathbf{F}_2 \mathbf{L}_1 \mathbf{F}_2^{-1} + \mathbf{F}_2^{-T} \mathbf{L}_1^T \mathbf{F}_2^T\right)\right] - \mu_p \mathbf{B}_p.(\mathbf{D} - \mathbf{D}_1) = \eta_1' \mathbf{D}_1.\mathbf{D}_1 + \eta_2' \mathbf{D}_2.\mathbf{D}_2, \tag{3.20}$$

which on further simplification leads to

$$\begin{aligned}(\mathbf{T} - \mu_3 \mathbf{B}_3 - \mu_p \mathbf{B}_p).\mathbf{D} + \mu_3 \mathbf{B}_3.\mathbf{D}_2 + \mu_p \mathbf{B}_p.\mathbf{D}_1 + \frac{\mu_3}{2}\mathbf{B}_3.\left(\mathbf{F}_2 \mathbf{L}_1 \mathbf{F}_2^{-1} + \mathbf{F}_2^{-T} \mathbf{L}_1^T \mathbf{F}_2^T\right) \\ = \eta_1' \mathbf{D}_1.\mathbf{D}_1 + \eta_2' \mathbf{D}_2.\mathbf{D}_2.\end{aligned} \tag{3.21}$$

---

[3]One can also choose the rate of dissipation function to depend on the stretch i.e., of the form $\xi \equiv \xi(\mathbf{D}_1, \mathbf{D}_2, \mathbf{B}_3, \mathbf{B}_p)$. The resulting constitutive relations will be a variant of the relations obtained when $\xi$ is of the form given in eqn. (3.14). The constitutive relations obtained by using $\xi(\mathbf{D}_1, \mathbf{D}_2, \mathbf{B}_3, \mathbf{B}_p)$ have relaxation times which depend on the stretch. Upon linearization, the two constitutive relations take the same form. Rajagopal and Srinivasa have discussed this issue for the Maxwell fluid in [16].

We shall assume that the body can undergo only isochoric motions and so

(3.22) $$tr(\mathbf{D}) = 0.$$

Also, since the body can actually attain the two natural configurations, the incompressibility constraint implies that

(3.23) $$tr(\mathbf{D}_1) = 0, \quad tr(\mathbf{D}_2) = 0,$$

where $tr\,(.)$ is the trace of second order tensor.

Since the right hand side of eqn. (3.21) does not depend on $\mathbf{D}$, along with eqn. (3.22), we have

(3.24) $$\mathbf{T} = -p\mathbf{I} + \mu_3 \mathbf{B}_3 + \mu_p \mathbf{B}_p,$$

where $-p\mathbf{I}$ is the reaction stress due to the constraint of incompressibility. Hence, eqn. (3.21) reduces to

(3.25) $$\mu_3 \mathbf{B}_3.\mathbf{D}_2 + \mu_p \mathbf{B}_p.\mathbf{D}_1 + \frac{\mu_3}{2}\mathbf{B}_3.\left(\mathbf{F}_2\mathbf{L}_1\mathbf{F}_2^{-1} + \mathbf{F}_2^{-T}\mathbf{L}_1^T\mathbf{F}_2^T\right) = \eta_1'\mathbf{D}_1.\mathbf{D}_1 + \eta_2'\mathbf{D}_2.\mathbf{D}_2.$$

Following Rajagopal and Srinivasa [16], we maximize the rate of dissipation in eqn. (3.16) along with the constraints in eqn. (3.23), (3.25), by varying $\mathbf{D}_1, \mathbf{D}_2$ for fixed $\mathbf{B}_2, \mathbf{B}_3$. We maximize the auxillary function $\Phi$ defined by

(3.26)
$$\Phi := \eta_1'\mathbf{D}_1.\mathbf{D}_1 + \eta_2'\mathbf{D}_2.\mathbf{D}_2$$
$$+ \lambda_1 \left[\eta_1'\mathbf{D}_1.\mathbf{D}_1 + \eta_2'\mathbf{D}_2.\mathbf{D}_2 - \mu_3\mathbf{B}_3.\mathbf{D}_2 - \mu_p\mathbf{B}_p.\mathbf{D}_1 - \frac{\mu_3}{2}\mathbf{B}_3.\left(\mathbf{F}_2\mathbf{L}_1\mathbf{F}_2^{-1} + \mathbf{F}_2^{-T}\mathbf{L}_1^T\mathbf{F}_2^T\right)\right]$$
$$+ \lambda_2 \mathbf{I}.\mathbf{D}_1 + \lambda_3 \mathbf{I}.\mathbf{D}_2$$

Now, setting $\partial\Phi/\partial\mathbf{D}_2 = 0, \partial\Phi/\partial\mathbf{D}_1 = 0$, and dividing the resulting equations by $\lambda_1$ and $\lambda_2$ respectively, for $\lambda_1, \lambda_2 \neq 0$, we get (also see appendix)

(3.27)
$$\mu_3\mathbf{B}_3 = \left(\frac{\lambda_1 + 1}{\lambda_1}\right) 2\eta_2'\mathbf{D}_2 + \frac{\lambda_3}{\lambda_1}\mathbf{I},$$
$$\mu_p\mathbf{B}_p + \frac{\mu_3}{2}\left(\mathbf{F}_2^T\mathbf{B}_3\mathbf{F}_2^{-T} + \mathbf{F}_2^{-1}\mathbf{B}_3\mathbf{F}_2\right) = \left(\frac{\lambda_1 + 1}{\lambda_1}\right) 2\eta_1'\mathbf{D}_1 + \frac{\lambda_2}{\lambda_1}\mathbf{I}.$$

Using eqn. (3.27) in eqn. (3.25), we get

(3.28) $$\frac{\lambda_1 + 1}{\lambda_1} = \frac{1}{2} - \frac{\mu_3 \mathbf{B}_3.\mathbf{F}_2\mathbf{W}_1\mathbf{F}_2^{-1}}{2\eta_1'\mathbf{D}_1.\mathbf{D}_1 + 2\eta_2'\mathbf{D}_2.\mathbf{D}_2},$$

where $\mathbf{W}_1 := \frac{1}{2}\left(\mathbf{L}_1 - \mathbf{L}_1^T\right)$. Hence,

(3.29)
$$\mathbf{T} = -p\mathbf{I} + \mu_3\mathbf{B}_3 + \mu_p\mathbf{B}_p,$$
$$\frac{\mu_3}{2}\left(\mathbf{F}_2^T\mathbf{B}_3\mathbf{F}_2^{-T} + \mathbf{F}_2^{-1}\mathbf{B}_3\mathbf{F}_2\right) + \mu_p\mathbf{B}_p = -p'\mathbf{I} + \eta_1\mathbf{D}_1,$$
$$\mu_3\mathbf{B}_3 = -p''\mathbf{I} + \eta_2\mathbf{D}_2,$$

where $p'$, $p''$ are the Lagrange multipliers with

(3.30)
$$-p' = \frac{1}{3}\left[\mu_3 tr(\mathbf{B}_3) + \mu_p tr(\mathbf{B}_p)\right], \quad -p'' = \frac{1}{3}\mu_3 tr(\mathbf{B}_3),$$
$$\eta_1 = 2\left(\frac{\lambda_1 + 1}{\lambda_1}\right)\eta_1', \quad \eta_2 = 2\left(\frac{\lambda_1 + 1}{\lambda_1}\right)\eta_2'.$$

Now, eqns. (3.10), (3.11) can be re-written as

$$\text{(3.31)} \qquad \mathbf{D}_1 = -\frac{1}{2}\mathbf{V}_p^{-1} \overset{\nabla}{\mathbf{B}}_p \mathbf{V}_p^{-1}, \quad \mathbf{D}_2 = -\frac{1}{2}\mathbf{V}_3^{-1} \overset{\nabla_p}{\mathbf{B}}_3 \mathbf{V}_3^{-1}.$$

Using eqn. $(3.31)_b$ in eqn. $(3.29)_c$, and post-multiplying and pre-multiplying with $\mathbf{V}_3$, we have

$$\text{(3.32)} \qquad \mu_3 \mathbf{B}_3^2 = \frac{1}{3}\mu_3 tr(\mathbf{B}_3)\mathbf{B}_3 - \frac{\eta_2}{2} \overset{\nabla_p}{\mathbf{B}}_3.$$

In addition, using eqn. $(3.31)_a$ in eqn. $(3.29)_b$, post-multiplying and pre-multiplying with $\mathbf{V}_p$, and using eqn. (3.4), we get

$$\text{(3.33)} \qquad \frac{\mu_3}{2}(\mathbf{B}_p \mathbf{B}_3 + \mathbf{B}_3 \mathbf{B}_p) + \mu_p \mathbf{B}_p^2 = \frac{1}{3}[\mu_3 tr(\mathbf{B}_3) + \mu_p tr(\mathbf{B}_p)]\mathbf{B}_p - \frac{\eta_1}{2} \overset{\nabla}{\mathbf{B}}_p.$$

Notice, from eqns. (3.32), (3.33), that the evolution of the natural configurations $\kappa_{p_1(t)}$ and $\kappa_{p_2(t)}$ are coupled. These two equations are to be solved simultaneously to determine their evolution. We shall denote $\mu_3 \mathbf{B}_3, \mu_p \mathbf{B}_p$ by $\mathbf{S}_1, \mathbf{S}_2$ respectively. Then, the final constitutive relations – eqns. $(3.29)_a$, (3.32), (3.33) – reduce to

$$\text{(3.34)} \qquad \begin{aligned} \mathbf{T} &= -p\mathbf{I} + \mathbf{S}_1 + \mathbf{S}_2, \\ \mathbf{S}_1^2 &= \frac{1}{3}tr(\mathbf{S}_1)\mathbf{S}_1 - \frac{\eta_2}{2} \overset{\nabla_p}{\mathbf{S}}_1, \\ \frac{1}{2}(\mathbf{S}_2 \mathbf{S}_1 + \mathbf{S}_1 \mathbf{S}_2) + \mathbf{S}_2^2 &= \frac{1}{3}[tr(\mathbf{S}_1) + tr(\mathbf{S}_2)]\mathbf{S}_2 - \frac{\eta_1}{2} \overset{\nabla}{\mathbf{S}}_2. \end{aligned}$$

In the next sub-section, we shall show that the above constitutive model reduces to Burgers' model in one dimension.

3.3. **Reduction of the model to one dimensional Burgers' model.** In this sub-section, we shall first linearize the constitutive model given by eqn. (3.29) (we shall use eqn. (3.29), here, instead of eqn. (3.34), for the sake of simplicity) by assuming the elastic response is small (we shall define what we mean by small, precisely, later). Then we shall show that, in one dimension, the equations reduce to the one dimensional linear model due to Burgers (see eqn. (1.1)).

Now, eqn. $(3.29)_c$ can be re-written as

$$\text{(3.35)} \qquad \mu_3(\mathbf{B}_3 - \mathbf{I}) = \mu_3\left[\frac{1}{3}tr(\mathbf{B}_3) - 1\right]\mathbf{I} + \eta_2 \mathbf{D}_2.$$

If the displacement gradient with elastic response is small, i.e.,

$$\text{(3.36)} \qquad \max_{\substack{\mathbf{X} \in \mathcal{B} \\ t \in \mathcal{R}}} \left\| \frac{\partial \mathbf{u}(\mathbf{X}, t)}{\partial \mathbf{X}} \right\| = \mathbf{O}(\gamma), \quad \gamma \ll 1,$$

then

$$\text{(3.37)} \qquad \|\mathbf{B}_i - \mathbf{I}\| = \mathbf{O}(\gamma), \quad \gamma \ll 1, \quad i = 3, p,$$

and hence

$$\text{(3.38)} \qquad tr(\mathbf{B}_i) = 3 + \mathbf{O}(\gamma^2), \quad i = 3, p,$$

and so the first term on the right hand side of eqn. (3.35) can be dropped for small strain and eqn. (3.35) reduces to

$$\text{(3.39)} \qquad \mu_3(\mathbf{B}_3 - \mathbf{I}) = \eta_2 \mathbf{D}_2.$$

If $\lambda_i$ ($i = 1, 2, 3, p$ or no subscript) is the stretch, in one dimension, corresponding to the deformation gradient $\mathbf{F}_i$, then $\lambda_i^2$ and $\dot{\lambda}_i/\lambda_i$ are the equivalent values in one dimension, corresponding to $\mathbf{B}_i$ and $\mathbf{D}_i$. If $\epsilon_i$ is the true strain for the stretch $\lambda_i$, then $\epsilon_i = ln\lambda_i$ and so, $\dot{\epsilon}_i = \dot{\lambda}_i/\lambda_i$. Hence, eqn. (3.39) reduces to

$$\mu_3 \left(\lambda_3^2 - 1\right) = \eta_2 \frac{\dot{\lambda}_2}{\lambda_2} \tag{3.40}$$

or

$$\mu_3 \left(e^{2\epsilon_3} - 1\right) = \eta_2 \dot{\epsilon}_2, \tag{3.41}$$

which under the assumption of small strain (i.e., $\epsilon_3 \ll 1$), reduces to

$$2\mu_3 \epsilon_3 = \eta_2 \dot{\epsilon}_2. \tag{3.42}$$

Following a similar analysis, in one dimension, eqn. $(3.29)_b$ becomes

$$2\mu_3 \epsilon_3 + 2\mu_p \epsilon_p = \eta_1 \dot{\epsilon}_1, \tag{3.43}$$

and eqn. $(3.29)_a$ reduces to

$$\sigma = 2\mu_3 \epsilon_3 + 2\mu_p \epsilon_p, \tag{3.44}$$

where $\sigma$ is the one dimensional stress. In addition, eqn. (3.3), reduces to

$$\lambda = \lambda_1 \lambda_2 \lambda_3, \tag{3.45}$$

and so

$$\epsilon = \epsilon_1 + \epsilon_2 + \epsilon_3. \tag{3.46}$$

Similarly, eqns. (3.4), (3.5) reduce to

$$\epsilon = \epsilon_p + \epsilon_1, \quad \epsilon_p = \epsilon_2 + \epsilon_3. \tag{3.47}$$

The eqns. (3.42–3.44), (3.46), (3.47) are, in fact, the equations obtained if we have the spring-dashpot arrangement shown in Fig. 3(a).

We shall now show that these equations (i.e., eqns. (3.42–3.44), (3.46), (3.47)) reduce to the form of eqn. (1.1). Now, differentiating eqn. $(3.47)_b$ with respect to time and using eqn. (3.42), we obtain

$$\dot{\epsilon}_p = \frac{2\mu_3}{\eta_2} \epsilon_3 + \dot{\epsilon}_3. \tag{3.48}$$

Also, differentiating eqn. (3.44) with respect to time and dividing by $\eta_2$, we find

$$\frac{\dot{\sigma}}{\eta_2} = \frac{2\mu_3}{\eta_2} \dot{\epsilon}_3 + \frac{2\mu_p}{\eta_2} \dot{\epsilon}_p, \tag{3.49}$$

and differentiating eqn. (3.44) twice with respect to time and dividing by $2\mu_3$ leads to

$$\frac{\ddot{\sigma}}{2\mu_3} = \ddot{\epsilon}_3 + \frac{\mu_p}{\mu_3} \ddot{\epsilon}_p. \tag{3.50}$$

We add eqns. (3.49), (3.50) and use eqn. (3.48), to get

$$\frac{\dot{\sigma}}{\eta_2} + \frac{\ddot{\sigma}}{2\mu_3} = \frac{2\mu_p}{\eta_2} \dot{\epsilon}_p + \left(1 + \frac{\mu_p}{\mu_3}\right) \ddot{\epsilon}_p. \tag{3.51}$$

From eqns. (3.43), (3.44)

(3.52) $$\sigma = \eta_1 \dot{\epsilon}_1 = \eta_1 (\dot{\epsilon} - \dot{\epsilon}_p), \quad \text{or,} \quad \dot{\epsilon}_p = \dot{\epsilon} - \frac{\sigma}{\eta_1}.$$

Using eqn. (3.52) in eqn. (3.51) leads to

(3.53) $$\frac{2\mu_p}{\eta_1 \eta_2} \sigma + \left( \frac{1}{\eta_1} + \frac{\mu_p}{\mu_3 \eta_1} + \frac{1}{\eta_2} \right) \dot{\sigma} + \frac{\ddot{\sigma}}{2\mu_3} = \frac{2\mu_p}{\eta_2} \dot{\epsilon} + \left( 1 + \frac{\mu_p}{\mu_3} \right) \ddot{\epsilon},$$

which can be re-written as

(3.54) $$\sigma + \left( \frac{\eta_2}{2\mu_p} + \frac{\eta_2}{2\mu_3} + \frac{\eta_1}{2\mu_p} \right) \dot{\sigma} + \frac{\eta_1 \eta_2}{4\mu_p \mu_3} \ddot{\sigma} = \eta_1 \dot{\epsilon} + \frac{\eta_1 \eta_2}{2\mu_p} \left( 1 + \frac{\mu_p}{\mu_3} \right) \ddot{\epsilon}.$$

Eqn. (3.54) is in the same form as eqn. (1.1), with

(3.55) $$p_1 = \frac{\eta_2}{2\mu_p} + \frac{\eta_2}{2\mu_3} + \frac{\eta_1}{2\mu_p}, \quad p_2 = \frac{\eta_1 \eta_2}{4\mu_p \mu_3}, \quad q_1 = \eta_1, \quad q_2 = \frac{\eta_1 \eta_2}{2\mu_p} \left( 1 + \frac{\mu_p}{\mu_3} \right).$$

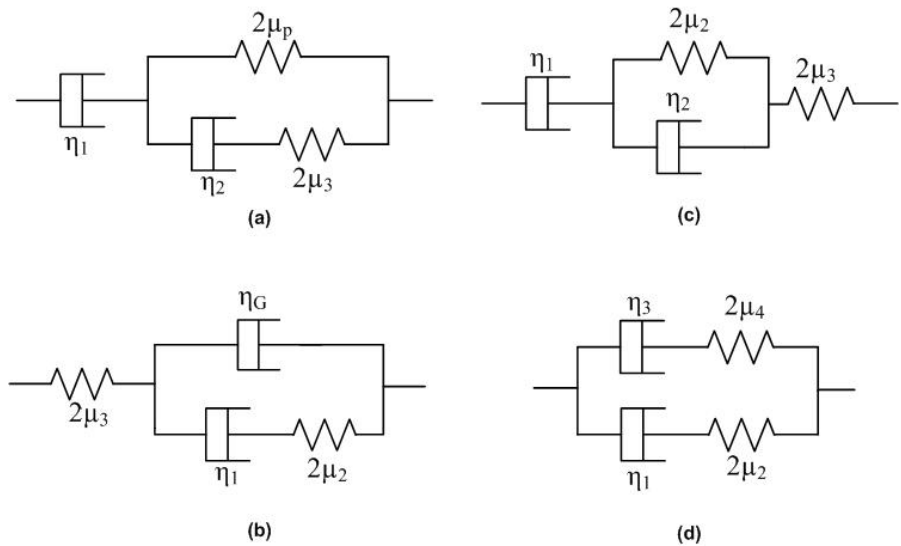

FIGURE 3. Various spring-dashpot arrangements which reduce to the one-dimensional Burgers' fluid model (eqn. (1.1)).

## 4. MODEL 2

4.1. **Preliminaries.** Once again, let $\kappa_R$ denote the reference configuration of the body. We shall assume that the body has two evolving natural configurations (denoted by $\kappa_{p_1(t)}, \kappa_{p_2(t)}$), but the manner in which they store the energy is different from that considered previously, with $\mathbf{F}_i, i = 1, 2, 3$, being the gradients of the motion as discussed in model 1. We shall also use the definitions in eqns. (3.1), (3.2). Thus, eqn. (3.3) applies here too. In addition, let us call the gradient of the motion from $\kappa_R$ to $\kappa_{p_2(t)}$ by $\mathbf{F}_G$ (see Fig. 4). It immediately follows that

(4.1) $$\mathbf{F}_G = \mathbf{F}_2 \mathbf{F}_1.$$

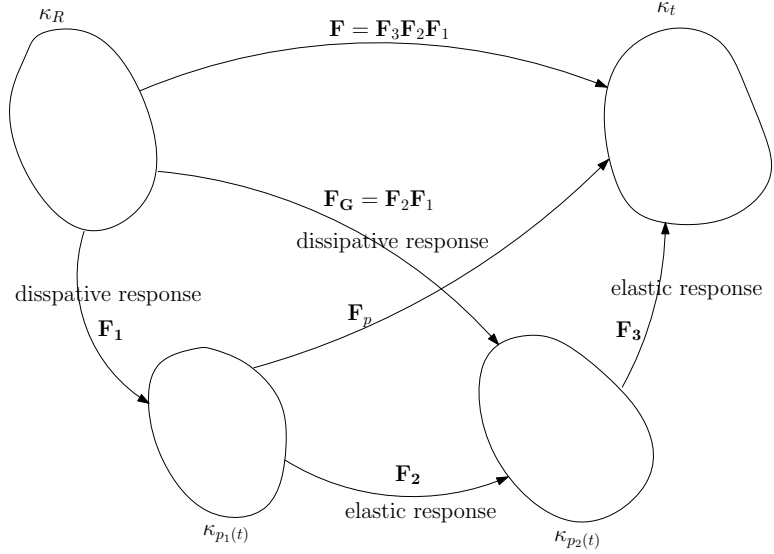

FIGURE 4. Schematic to illustrate the natural configurations for model 2. The body dissipates like a viscous fluid during its motion from, $\kappa_R$ to $\kappa_{p_2(t)}$, and $\kappa_R$ to $\kappa_{p_1(t)}$. The body stores energy like a neo-Hookean solid during its motion from $\kappa_{p_1(t)}$ to $\kappa_{p_2(t)}$ and $\kappa_{p_2(t)}$ to $\kappa_t$.

We shall denote the velocity gradient and its symmetric part corresponding to $\mathbf{F}_G$ by

$$(4.2) \qquad \mathbf{L}_G := \dot{\mathbf{F}}_G \mathbf{F}_G^{-1}, \quad \mathbf{D}_G := \frac{1}{2}\left(\mathbf{L}_G + \mathbf{L}_G^T\right).$$

Also, from eqns. (3.3), (4.1),

$$(4.3) \qquad \mathbf{F} = \mathbf{F}_3 \mathbf{F}_G.$$

Following a procedure similar to the one followed previously for model 1, it can be shown that

$$(4.4) \qquad \mathbf{D}_G = \mathbf{D}_2 + \frac{1}{2}(\mathbf{F}_2 \mathbf{L}_1 \mathbf{F}_2^{-1} + \mathbf{F}_2^{-T}\mathbf{L}_1^T \mathbf{F}_2^T), \quad \mathbf{D} = \mathbf{D}_3 + \frac{1}{2}(\mathbf{F}_3 \mathbf{L}_G \mathbf{F}_3^{-1} + \mathbf{F}_3^{-T}\mathbf{L}_G^T \mathbf{F}_3^T),$$

along with

$$(4.5) \qquad \overset{\triangledown}{\mathbf{B}}_3 = -2\mathbf{F}_3 \mathbf{D}_G \mathbf{F}_3^T, \quad \overset{\triangledown_G}{\mathbf{B}}_2 = -2\mathbf{F}_2 \mathbf{D}_1 \mathbf{F}_2^T,$$

where $\overset{\triangledown_G}{\mathbf{A}} := \dot{\mathbf{A}} - \mathbf{A}\mathbf{L}_G^T - \mathbf{L}_G \mathbf{A}$ is the Oldroyd derivative when the natural configuration $\kappa_{p_2(t)}$ is the current configuration. In addition, from eqns. (4.4) and (4.5), along with the assumption that $\mathbf{F}_2 = \mathbf{V}_2$, $\mathbf{F}_3 = \mathbf{V}_3$ in virtue of the body being isotropic, we get

$$(4.6) \qquad \mathbf{I}.\dot{\mathbf{B}}_2 = 2\mathbf{B}_2.(\mathbf{D}_G - \mathbf{D}_1), \quad \mathbf{I}.\dot{\mathbf{B}}_3 = 2\mathbf{B}_3.(\mathbf{D} - \mathbf{D}_G).$$

These relations should suffice for our calculations for studying the response of model 2.

4.2. **Constitutive assumptions.** In this model, we shall assume $\psi$, and $\xi$, to be of the form

$$(4.7) \qquad \psi \equiv \psi(\mathbf{B}_2, \mathbf{B}_3), \quad \xi \equiv \xi(\mathbf{D}_1, \mathbf{D}_G).$$

Now, assuming that the instantaneous elastic responses are isotropic and the body is incompressible, we choose

$$\psi(\mathbf{B}_2, \mathbf{B}_3) = \frac{\mu_2}{2\rho}(\mathbf{I}.\mathbf{B}_2 - 3) + \frac{\mu_3}{2\rho}(\mathbf{I}.\mathbf{B}_3 - 3), \tag{4.8}$$

and

$$\xi(\mathbf{D}_1, \mathbf{D}_G) = \eta_1 \mathbf{D}_1.\mathbf{D}_1 + \eta_G \mathbf{D}_G.\mathbf{D}_G. \tag{4.9}$$

The above assumption implies that the body possesses instantaneous elastic response from the current configuration $\kappa_t$ to the natural configuration $\kappa_{p_2(t)}$ and from the natural configuration $\kappa_{p_1(t)}$ to the other natural configuration $\kappa_{p_2(t)}$. It stores energy like a neo-Hookean solid during these two motions. In addition, the responses from the two natural configurations $(\kappa_{p_1(t)}, \kappa_{p_1(t)})$ to the reference configuration $\kappa_R$ are purely dissipative, similar to a linear viscous fluid. In fact, the response of the body as it moves from $\kappa_R$ to $\kappa_{p_2(t)}$ is similar to that of a "variant" of an Oldroyd-B fluid (see Rajagopal and Srinivasa [16]) i.e., the natural configuration $\kappa_{p_2(t)}$ evolves like that of an Oldroyd-B fluid with respect to the reference configuration $\kappa_R$.

On substituting eqns. (4.8), (4.9) in (2.12), using eqn. (4.6) and simplifying, we get

$$(\mathbf{T} - \mu_3 \mathbf{B}_3).\mathbf{D} + (\mu_3 \mathbf{B}_3 - \mu_2 \mathbf{B}_2).\mathbf{D}_G + \mu_2 \mathbf{B}_2.\mathbf{D}_1 = \eta_1 \mathbf{D}_1.\mathbf{D}_1 + \eta_G \mathbf{D}_G.\mathbf{D}_G. \tag{4.10}$$

Since, the right hand side of eqn. (4.10) does not depend on $\mathbf{D}$, the incompressibility constraint, $tr(\mathbf{D}) = 0$, leads to

$$\mathbf{T} = -p\mathbf{I} + \mu_3 \mathbf{B}_3, \tag{4.11}$$

where $-p\mathbf{I}$ is the reaction stress due to the incompressibility constraint. Using, eqn. (4.11) in (4.10), we must have

$$(\mu_3 \mathbf{B}_3 - \mu_2 \mathbf{B}_2).\mathbf{D}_G + \mu_2 \mathbf{B}_2.\mathbf{D}_1 = \eta_1 \mathbf{D}_1.\mathbf{D}_1 + \eta_G \mathbf{D}_G.\mathbf{D}_G. \tag{4.12}$$

Now, we maximize the rate of dissipation by varying $\mathbf{D}_1, \mathbf{D}_G$ for fixed $\mathbf{B}_2, \mathbf{B}_3$ with the constraints

$$tr(\mathbf{D}_1) = 0, \quad tr(\mathbf{D}_G) = 0. \tag{4.13}$$

Finally, we arrive at the following set of equations:

$$\begin{aligned} \mathbf{T} &= -p\mathbf{I} + \mu_3 \mathbf{B}_3, \\ \mu_3 \mathbf{B}_3 - \mu_2 \mathbf{B}_2 &= -p'\mathbf{I} + \eta_G \mathbf{D}_G, \\ \mu_2 \mathbf{B}_2 &= -p''\mathbf{I} + \eta_1 \mathbf{D}_1, \end{aligned} \tag{4.14}$$

where $p$, $p'$, $p''$ are the Lagrange multipliers with

$$\begin{aligned} p' &= -\frac{1}{3}\left[\mu_3 tr(\mathbf{B}_3) - \mu_2 tr(\mathbf{B}_2)\right], \\ p'' &= -\frac{1}{3}\mu_2 tr(\mathbf{B}_2). \end{aligned} \tag{4.15}$$

Pre-multiplying and post-multiplying eqns. (4.14)$_b$, (4.14)$_c$ by $\mathbf{V}_3$ and $\mathbf{V}_2$ respectively, eqn. (4.14) reduces to

(4.16)
$$\mathbf{T} = -p\mathbf{I} + \mu_3 \mathbf{B}_3,$$
$$\mu_3 \mathbf{B}_3^2 - \mu_2 \mathbf{V}_3 \mathbf{B}_2 \mathbf{V}_3 = -p' \mathbf{B}_3 - \frac{\eta_G}{2} \overset{\nabla}{\mathbf{B}}_3,$$
$$\mu_2 \mathbf{B}_2^2 = -p'' \mathbf{B}_2 - \frac{\eta_1}{2} \overset{\nabla_G}{\mathbf{B}}_2,$$

with eqn. (4.15). If we denote $\mu_2 \mathbf{B}_2, \mu_3 \mathbf{B}_3$ by $\mathbf{S}_1, \mathbf{S}_2$ respectively, then the final constitutive relations for this model are

(4.17)
$$\mathbf{T} = -p\mathbf{I} + \mathbf{S}_2,$$
$$\mathbf{S}_1^2 = \frac{1}{3} tr(\mathbf{S}_1) \mathbf{S}_1 - \frac{\eta_1}{2} \overset{\nabla_G}{\mathbf{S}}_1,$$
$$\mathbf{S}_2^2 - \sqrt{\mathbf{S}_2} \mathbf{S}_1 \sqrt{\mathbf{S}_2} = \frac{1}{3} [tr(\mathbf{S}_2) - tr(\mathbf{S}_1)] \mathbf{S}_2 - \frac{\eta_G}{2} \overset{\nabla}{\mathbf{S}}_2.$$

**4.3. Reduction of the model to one dimensional Burgers' model.** For simplicity, we shall use eqn. (4.14) for the reduction. Now, eqns. (4.14)$_{b,c}$ can be re-written as

(4.18)
$$\mu_3 (\mathbf{B}_3 - \mathbf{I}) - \mu_2 (\mathbf{B}_2 - \mathbf{I}) = \frac{1}{3} [\mu_3 (tr(\mathbf{B}_3) - 3) - \mu_2 (tr(\mathbf{B}_2) - 3)] \mathbf{I} + \eta_G \mathbf{D}_G,$$
$$\mu_2 (\mathbf{B}_2 - \mathbf{I}) = \frac{1}{3} \mu_2 (tr(\mathbf{B}_2) - 3) \mathbf{I} + \eta_1 \mathbf{D}_1.$$

Assuming that the displacement gradient associated with elastic response is small, leads to

(4.19)
$$\|\mathbf{B}_i - \mathbf{I}\| = \mathbf{O}(\gamma), \quad \gamma \ll 1, \quad i = 2, 3.$$

The first term on the right hand sides of eqn. (4.18)$_{b,c}$ can be neglected. Then, eqn. (4.18) reduces to

(4.20)
$$\mu_3 (\mathbf{B}_3 - \mathbf{I}) - \mu_2 (\mathbf{B}_2 - \mathbf{I}) = \eta_G \mathbf{D}_G,$$
$$\mu_2 (\mathbf{B}_2 - \mathbf{I}) = \eta_1 \mathbf{D}_1.$$

In one dimension, eqn. (4.20) becomes

(4.21)
$$\mu_3(\lambda_3^2 - 1) - \mu_2(\lambda_2^2 - 1) = \eta_G \frac{\dot{\lambda}_G}{\lambda_G}, \quad \mu_2(\lambda_2^2 - 1) = \eta_1 \frac{\dot{\lambda}_1}{\lambda_1},$$

where $\lambda_i$ ($i = 1, 2, 3, G$ or no subscript) is the stretch, in one dimension, corresponding to the right stretch tensor $\mathbf{V}_i$. Using $ln\lambda_i = \epsilon_i$ ($\epsilon$ is the true strain), with the assumption of $\epsilon_i \ll 1$, eqn. (4.21) reduces to

(4.22)
$$2\mu_3 \epsilon_3 - 2\mu_2 \epsilon_2 = \eta_G \dot{\epsilon}_G, \quad 2\mu_2 \epsilon_2 = \eta_1 \dot{\epsilon}_1.$$

In addition, eqn. (4.14)$_a$ reduces to

(4.23)
$$\sigma = 2\mu_3 \epsilon_3.$$

Eqn. (4.3) together with eqn. (4.1), in one dimension, reduces to

(4.24)
$$\epsilon = \epsilon_G + \epsilon_3, \quad \text{or}, \quad \epsilon_G = \epsilon_2 + \epsilon_1.$$

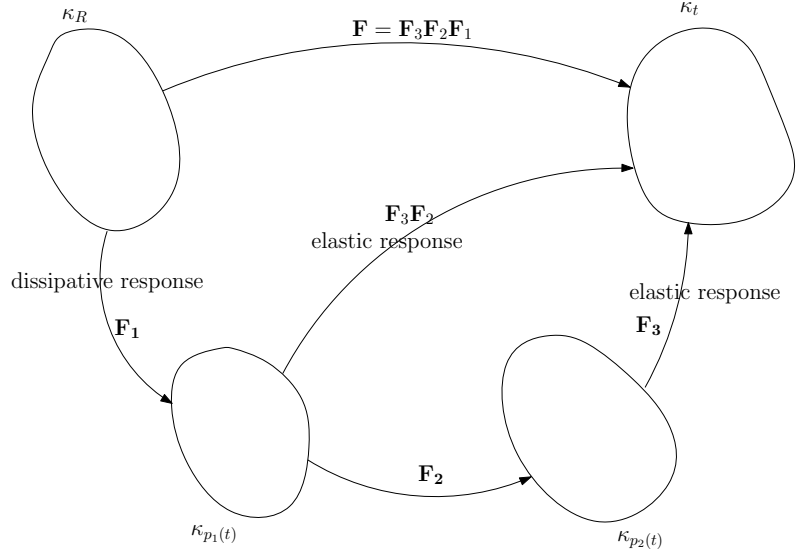

FIGURE 5. Schematic to illustrate the natural configurations for model 3. The body's response is viscous fluid-like and elastic solid-like, during its motion from, $\kappa_R$ to $\kappa_{p_t(t)}$, and $\kappa_{p_2(t)}$ to $\kappa_t$ respectively. From $\kappa_{p_1(t)}$ to $\kappa_{p_2(t)}$, the response is Kelvin-Voigt solid-like.

The spring-dashpot arrangement in Fig. 3(b) also yields eqns. (4.22), (4.23) along with eqns. (3.46), (4.24). These equations on simplification reduce to

$$\sigma + \left(\frac{\eta_1}{2\mu_2} + \frac{\eta_1}{2\mu_3} + \frac{\eta_G}{2\mu_3}\right)\dot{\sigma} + \frac{\eta_1\eta_G}{4\mu_2\mu_3}\ddot{\sigma} = (\eta_1 + \eta_G)\dot{\epsilon} + \frac{\eta_1\eta_G}{2\mu_2}\ddot{\epsilon}, \tag{4.25}$$

which is same as the Burgers' one dimensional model (eqn. (1.1)), with

$$p_1 = \frac{\eta_1}{2\mu_2} + \frac{\eta_1}{2\mu_3} + \frac{\eta_G}{2\mu_3}, \quad p_2 = \frac{\eta_1\eta_G}{4\mu_2\mu_3}, \quad q_1 = \eta_1 + \eta_G, \quad q_2 = \frac{\eta_1\eta_G}{2\mu_2}. \tag{4.26}$$

## 5. MODEL 3

5.1. **Preliminaries.** As with models 1 and 2, for this model we shall assume that the body has two evolving natural configurations ($\kappa_{p_1(t)}, \kappa_{p_1(t)}$, see Fig. 5). We shall also use the definition of $\mathbf{F}_i, i = 1, 2, 3$ used for the previous models in addition to the definitions in eqns. (3.1), (3.2) and the relation eqn. (3.3). Further, we shall also choose $\mathbf{F}_2 = \mathbf{V}_2$ and $\mathbf{F}_3 = \mathbf{V}_3$. We recall from the preliminary discussion concerning models 1 and 2, that

$$\mathbf{D}_2 = -\frac{1}{2}\mathbf{V}_3^{-1}\overset{\nabla_p}{\mathbf{B}}_3\mathbf{V}_3^{-1}, \quad \mathbf{D}_1 = -\frac{1}{2}\mathbf{V}_2^{-1}\overset{\nabla_G}{\mathbf{B}}_2\mathbf{V}_2^{-1}. \tag{5.1}$$

These definitions and relations shall be used in the following analysis.

5.2. **Constitutive assumptions.** For this model, we shall assume the specific stored energy, $\psi$ and the rate of dissipation, $\xi$ to be of the form

$$\psi \equiv \psi(\mathbf{B}_2, \mathbf{B}_3), \quad \xi \equiv \xi(\mathbf{D}_1, \mathbf{D}_2). \tag{5.2}$$

Specifically, in virtue of the body being incompressible and isotropic, we choose,

$$\psi(\mathbf{B}_2, \mathbf{B}_3) = \frac{\mu_2}{2\rho}(\mathbf{I}.\mathbf{B}_2 - 3) + \frac{\mu_3}{2\rho}(\mathbf{I}.\mathbf{B}_3 - 3), \tag{5.3}$$

and

$$\xi(\mathbf{D}_1, \mathbf{D}_2) = \eta_1' \mathbf{D}_1.\mathbf{D}_1 + \eta_2' \mathbf{D}_2.\mathbf{D}_2, \tag{5.4}$$

i.e., the body possesses instantaneous elastic response from the current configuration $\kappa_t$ to the natural configuration $\kappa_{p_2(t)}$ and stores energy like a neo-Hookean solid. Also, the response of the body between $\kappa_{p_1(t)}$ to $\kappa_{p_2(t)}$ is similar to that of a Kelvin-Voigt solid. The body also dissipates like a linear viscous fluid during its motion from $\kappa_R$ to $\kappa_{p_1(t)}$.

On substituting eqn. (5.3) into eqn. (2.12) and using eqn. (3.18) we get,

$$\mathbf{T}.\mathbf{D} - \mu_2 \mathbf{B}_2.\mathbf{D}_2 - \mu_3 \mathbf{B}_3.\left[\mathbf{D} - \mathbf{D}_2 - \frac{1}{2}\left(\mathbf{F}_2 \mathbf{L}_1 \mathbf{F}_2^{-1} + \mathbf{F}_2^{-T} \mathbf{L}_1^T \mathbf{F}_2^T\right)\right] = \eta_1' \mathbf{D}_1.\mathbf{D}_1 + \eta_2' \mathbf{D}_2.\mathbf{D}_2, \tag{5.5}$$

which reduces to

$$(\mathbf{T} - \mu_3 \mathbf{B}_3).\mathbf{D} + (\mu_3 \mathbf{B}_3 - \mu_2 \mathbf{B}_2).\mathbf{D}_2 + \frac{\mu_3}{2} \mathbf{B}_3.\left(\mathbf{F}_2 \mathbf{L}_1 \mathbf{F}_2^{-1} + \mathbf{F}_2^{-T} \mathbf{L}_1^T \mathbf{F}_2^T\right) = \eta_1' \mathbf{D}_1.\mathbf{D}_1 + \eta_2' \mathbf{D}_2.\mathbf{D}_2. \tag{5.6}$$

Using eqn. (5.6), we maximize the rate of dissipation with incompressibility as a constraint, i.e.,

$$tr(\mathbf{D}) = tr(\mathbf{D}_1) = tr(\mathbf{D}_2) = 0, \tag{5.7}$$

by varying $\mathbf{D}, \mathbf{D}_1, \mathbf{D}_2$ for fixed $\mathbf{B}_2, \mathbf{B}_3$ and get:

$$\begin{aligned} \mathbf{T} &= -p\mathbf{I} + \mu_3 \mathbf{B}_3, \\ \mu_3 \mathbf{B}_3 - \mu_2 \mathbf{B}_2 &= -p'\mathbf{I} + \eta_2 \mathbf{D}_2, \\ \frac{\mu_3}{2}\left(\mathbf{F}_2^T \mathbf{B}_3 \mathbf{F}_2^{-T} + \mathbf{F}_2^{-1} \mathbf{B}_3 \mathbf{F}_2\right) &= -p''\mathbf{I} + \eta_1 \mathbf{D}_1, \end{aligned} \tag{5.8}$$

where $p$, $p'$, $p''$ are the Lagrange multipliers with

$$\begin{aligned} -p' &= \frac{1}{3}\left[\mu_3 tr(\mathbf{B}_3) - \mu_2 tr(\mathbf{B}_2)\right], \\ -p'' &= \frac{1}{3}\mu_3 tr(\mathbf{B}_3), \end{aligned} \tag{5.9}$$

and

$$\eta_i = \eta_i' \left(1 - \frac{\mu_3 \mathbf{B}_3.\mathbf{F}_2 \mathbf{W}_1 \mathbf{F}_2^{-1}}{\eta_1' \mathbf{D}_1.\mathbf{D}_1 + \eta_2' \mathbf{D}_2.\mathbf{D}_2}\right), \quad i = 1, 2. \tag{5.10}$$

Pre-multiplying and post-multiplying eqn. (5.8)$_b$ by $\mathbf{V}_3$, pre-multiplying and post-multiplying eqn. (5.8)$_c$ by $\mathbf{V}_2$ and using eqn. (5.1), we find that

$$\begin{aligned} \mathbf{T} &= -p\mathbf{I} + \mu_3 \mathbf{B}_3, \\ \mu_3 \mathbf{B}_3^2 - \mu_2 \mathbf{V}_3 \mathbf{B}_2 \mathbf{V}_3 &= -p'\mathbf{B}_3 - \frac{\eta_2}{2} \overset{\nabla_p}{\mathbf{B}_3}, \\ \frac{\mu_3}{2}\left(\mathbf{B}_2 \mathbf{B}_3 + \mathbf{B}_3 \mathbf{B}_2\right) &= -p''\mathbf{B}_2 - \frac{\eta_1}{2} \overset{\nabla_G}{\mathbf{B}_2}, \end{aligned} \tag{5.11}$$

along with eqn. (5.9).

If we call $\mu_3 \mathbf{B}_3, \mu_2 \mathbf{B}_2$ by $\mathbf{S}_1, \mathbf{S}_2$ respectively, then, the final form for the constitutive relation can be given as

$$\mathbf{T} = -p\mathbf{I} + \mathbf{S}_1,$$

(5.12)
$$\mathbf{S}_1^2 - \sqrt{\mathbf{S}_1}\mathbf{S}_2\sqrt{\mathbf{S}_1} = \frac{1}{3}\left[tr(\mathbf{S}_1) - tr(\mathbf{S}_2)\right]\mathbf{S}_1 - \frac{\eta_2}{2}\overset{\nabla_p}{\mathbf{S}_1},$$

$$\frac{1}{2}\left(\mathbf{S}_2\mathbf{S}_1 + \mathbf{S}_1\mathbf{S}_2\right) = \frac{1}{3}tr(\mathbf{S}_1)\mathbf{S}_2 - \frac{\eta_1}{2}\overset{\nabla_G}{\mathbf{S}_2}.$$

5.3. **Reduction of the model to one dimensional Burgers' model.** Following the method used in 4.3, eqn. (5.8), in one dimension, reduces to

(5.13)
$$\sigma = 2\mu_3 \epsilon_3,$$
$$2\mu_3 \epsilon_3 - 2\mu_2 \epsilon_2 = \eta_2 \dot{\epsilon}_2,$$
$$2\mu_3 \epsilon_3 = \eta_1 \dot{\epsilon}_1.$$

The above set of equations, can also be obtained from the spring-dashpot arrangement in Fig. 3(c).

Now, eqn. (5.13) can be re-written as

(5.14) $\qquad \sigma = 2\mu_3 \epsilon_3, \quad \sigma = 2\mu_2 \epsilon_2 + \eta_2 \dot{\epsilon}_2, \quad \sigma = \eta_1 \dot{\epsilon}_1.$

Also, differentiating eqn. (3.46) with respect to time and using eqns. (5.14)$_{a,c}$, we obtain

(5.15) $\qquad \dot{\epsilon} = \dfrac{\sigma}{\eta_1} + \dfrac{\dot{\sigma}}{2\mu_3} + \dot{\epsilon}_2.$

Now, multiplying eqn. (5.15) with $2\mu_2$, multiplying the derivative of eqn. (5.15) with respect to time with $\eta_2$; then, adding these two equations, along with eqn. (5.14)$_b$, we get

(5.16) $\qquad \dfrac{2\mu_2}{\eta_1}\sigma + \left(1 + \dfrac{\eta_2}{\eta_1} + \dfrac{\mu_2}{\mu_3}\right)\dot{\sigma} + \dfrac{\eta_2}{2\mu_3}\ddot{\sigma} = 2\mu_2 \dot{\epsilon} + \eta_2 \ddot{\epsilon},$

re-written as

(5.17) $\qquad \sigma + \left(\dfrac{\eta_1}{2\mu_2} + \dfrac{\eta_2}{2\mu_2} + \dfrac{\eta_1}{2\mu_3}\right)\dot{\sigma} + \dfrac{\eta_1 \eta_2}{4\mu_2 \mu_3}\ddot{\sigma} = \eta_1 \dot{\epsilon} + \dfrac{\eta_1 \eta_2}{2\mu_2}\ddot{\epsilon}.$

Thus, eqn. (5.17) has the same form as eqn. (1.1), with

(5.18) $\qquad p_1 = \dfrac{\eta_1}{2\mu_2} + \dfrac{\eta_2}{2\mu_2} + \dfrac{\eta_1}{2\mu_3}, \quad p_2 = \dfrac{\eta_1 \eta_2}{4\mu_2 \mu_3}, \quad q_1 = \eta_1, \quad q_2 = \dfrac{\eta_1 \eta_2}{2\mu_2}.$

6. MODEL 4

6.1. **Preliminaries.** Once again, we shall assume that the body has two natural configurations associated with it, denoted by $\kappa_{p_1(t)}, \kappa_{p_2(t)}$. However, in this model, the evolution equations of the two natural configurations are not coupled and they evolve independently (see Fig. 6). We shall denote the gradients of the motion from $\kappa_R$ to $\kappa_{p_1(t)}$ and from $\kappa_{p_1(t)}$ to $\kappa_t$ by $\mathbf{F}_1, \mathbf{F_2}$. We shall also denote the gradients of the motion from $\kappa_R$ to $\kappa_{p_2(t)}$ and from $\kappa_{p_2(t)}$ to $\kappa_t$ by $\mathbf{F}_3, \mathbf{F_4}$. It follows that

(6.1) $\qquad \mathbf{F} = \mathbf{F}_2 \mathbf{F}_1 = \mathbf{F}_4 \mathbf{F}_3.$

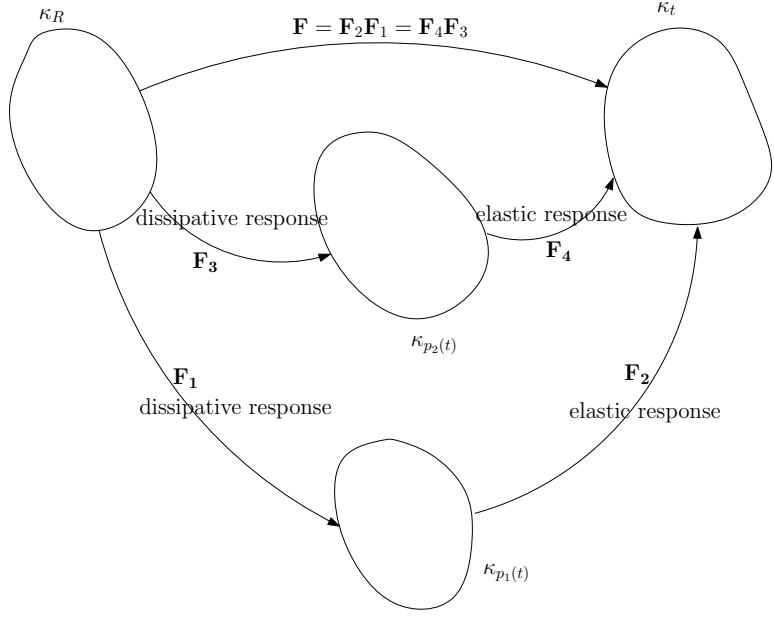

FIGURE 6. Schematic to illustrate the natural configurations for model 4. The body's response is similar to that of a "mixture" of two Maxwell-like fluids with different relaxation times.

The left stretch tensor, velocity gradient and its corresponding symmetric part are denoted by

$$\text{(6.2)} \qquad \mathbf{B}_i := \mathbf{F}_i \mathbf{F}_i^T, \quad \mathbf{L}_i := \dot{\mathbf{F}}_i \mathbf{F}_i^{-1}, \quad \mathbf{D}_i := \frac{1}{2}\left(\mathbf{L}_i + \mathbf{L}_i^T\right), \quad i = 1, 2, 3, 4.$$

Also, a straightforward calculation leads to

$$\text{(6.3)} \qquad \overset{\triangledown}{\mathbf{B}_2} = -2\mathbf{F}_2 \mathbf{D}_1 \mathbf{F}_2^T, \quad \overset{\triangledown}{\mathbf{B}_4} = -2\mathbf{F}_4 \mathbf{D}_3 \mathbf{F}_4^T.$$

6.2. **Constitutive assumptions.** Here, we shall assume the specific stored energy, $\psi$ and the rate of dissipation, $\xi$ to be of the form

$$\text{(6.4)} \qquad \psi \equiv \psi(\mathbf{B}_2, \mathbf{B}_4), \quad \xi \equiv \xi(\mathbf{D}_1, \mathbf{D}_3).$$

As the material is isotropic and incompressible, we choose,

$$\text{(6.5)} \qquad \psi(\mathbf{B}_2, \mathbf{B}_4) = \frac{\mu_2}{2\rho}(\mathbf{I}.\mathbf{B}_2 - 3) + \frac{\mu_4}{2\rho}(\mathbf{I}.\mathbf{B}_4 - 3),$$

and

$$\text{(6.6)} \qquad \xi(\mathbf{D}_1, \mathbf{D}_3) = \eta_1 \mathbf{D}_1.\mathbf{D}_1 + \eta_3 \mathbf{D}_3.\mathbf{D}_3.$$

This means that the response of the natural configurations $(\kappa_{p_1(t)}, \kappa_{p_2(t)})$ from the current configuration is like that of a neo-Hookean solid and the response from the reference configuration to the natural configurations is similar to that of a linear viscous fluid. Thus, Burgers' fluid can also be perceived as a "mixture" of two Maxwell-like fluids with different relaxation times.

We shall set

$$\text{(6.7)} \qquad \mathbf{F}_2 = \mathbf{V}_2, \quad \mathbf{F}_4 = \mathbf{V}_4,$$

where $\mathbf{V}_2, \mathbf{V}_4$ are the right stretch tensors in the polar decomposition of $\mathbf{F}_2, \mathbf{F}_4$, based on the assumption of isotropic elastic response. Hence, from eqn. (6.3) and eqn. (6.7), we have

(6.8) $$\mathbf{I}.\dot{\mathbf{B}}_2 = 2\mathbf{B}_2.(\mathbf{D}-\mathbf{D}_1), \quad \mathbf{I}.\dot{\mathbf{B}}_4 = 2\mathbf{B}_4.(\mathbf{D}-\mathbf{D}_3).$$

On entering eqns. (6.5), (6.6) in eqn. (2.12) and using eqn. (6.8), we get

(6.9) $$(\mathbf{T}-\mu_2\mathbf{B}_2-\mu_4\mathbf{B}_4).\mathbf{D} + \mu_2\mathbf{B}_2.\mathbf{D}_1 + \mu_4\mathbf{B}_4.\mathbf{D}_3 = \eta_1\mathbf{D}_1.\mathbf{D}_1 + \eta_3\mathbf{D}_3.\mathbf{D}_3.$$

Using the constraint of incompressibility

(6.10) $$tr(\mathbf{D}) = tr(\mathbf{D}_1) = tr(\mathbf{D}_3) = 0,$$

and eqn.(6.9), we maximize the rate of dissipation by varying $\mathbf{D}, \mathbf{D}_1, \mathbf{D}_3$ for fixed $\mathbf{B}_2, \mathbf{B}_4$ and get:

(6.11) $$\begin{aligned} \mathbf{T} &= -p\mathbf{I} + \mu_2\mathbf{B}_2 + \mu_4\mathbf{B}_4, \\ \mu_2\mathbf{B}_2 &= -p'\mathbf{I} + \eta_1\mathbf{D}_1, \\ \mu_4\mathbf{B}_4 &= -p''\mathbf{I} + \eta_3\mathbf{D}_3, \end{aligned}$$

where $p, p', p''$ are the Lagrange multipliers with

(6.12) $$-p' = \frac{1}{3}\mu_2 tr(\mathbf{B}_2), \quad -p'' = \frac{1}{3}\mu_4 tr(\mathbf{B}_4).$$

Pre-multiplying and post-multiplying, eqn. (6.11)$_b$ by $\mathbf{V}_2$, and eqn. (6.11)$_c$ by $\mathbf{V}_4$; then, using eqn. (6.3), we arrive at

(6.13) $$\begin{aligned} \mu_2\mathbf{B}_2^2 &= -p'\mathbf{B}_2 - \frac{\eta_1}{2}\stackrel{\triangledown}{\mathbf{B}}_2, \\ \mu_4\mathbf{B}_4^2 &= -p''\mathbf{B}_4 - \frac{\eta_3}{2}\stackrel{\triangledown}{\mathbf{B}}_4, \end{aligned}$$

Eqns. (6.13)$_{a,b}$ represents the evolution equations of the natural configurations $(\kappa_{p_1(t)}, \kappa_{p_2(t)}$ respectively). If we denote $\mu_2\mathbf{B}_2, \mu_4\mathbf{B}_4$ by $\mathbf{S}_1, \mathbf{S}_2$ respectively, then the final constitutive relations, for model 4, are

(6.14) $$\begin{aligned} \mathbf{T} &= -p\mathbf{I} + \mathbf{S}_1 + \mathbf{S}_2, \\ \mathbf{S}_1^2 &= \frac{1}{3}tr(\mathbf{S}_1)\mathbf{S}_1 - \frac{\eta_1}{2}\stackrel{\triangledown}{\mathbf{S}}_1, \\ \mathbf{S}_2^2 &= \frac{1}{3}tr(\mathbf{S}_2)\mathbf{S}_2 - \frac{\eta_1}{2}\stackrel{\triangledown}{\mathbf{S}}_2. \end{aligned}$$

This model is a variation of the model proposed by Murali Krishnan and Rajagopal [8]. They considered stretch dependent dissipation, in constrast to our linear viscous fluid type dissipation.

6.3. **Reduction of the model to the one dimensional Burgers' model.** For this model, we shall once again assume that the displacement gradient associated with the elastic response is small, and thus

(6.15) $$\|\mathbf{B}_i - \mathbf{I}\| = \mathbf{O}(\gamma), \quad \gamma \ll 1, \quad i = 2, 4.$$

Then, eqn. (6.11) becomes

(6.16) $$\begin{aligned} \mathbf{T} &= -p\mathbf{I} + \mu_2\mathbf{B}_2 + \mu_4\mathbf{B}_4, \\ \mu_2\left(\mathbf{B}_2 - \mathbf{I}\right) &= \eta_1\mathbf{D}_1, \\ \mu_4\left(\mathbf{B}_4 - \mathbf{I}\right) &= \eta_3\mathbf{D}_3, \end{aligned}$$

which in one dimension reduces to

(6.17)
$$\sigma = 2\mu_2 \epsilon_2 + 2\mu_4 \epsilon_4,$$
$$2\mu_2 \epsilon_2 = \eta_1 \dot{\epsilon}_1,$$
$$2\mu_4 \epsilon_4 = \eta_3 \dot{\epsilon}_3.$$

Further, eqn (6.1), in one dimension, reduces to

(6.18) $$\epsilon = \epsilon_2 + \epsilon_1 = \epsilon_3 + \epsilon_4.$$

In fact, the spring-dashpot arrangement Fig. 3(d) leads to eqn. (6.17), (6.18). We shall now show that these two equations, on simplification lead to eqn. (1.1). Differentiating eqn. (6.18) with respect to time and using eqn. $(6.17)_{b,c}$, we have

(6.19)
$$\dot{\epsilon} = \frac{2\mu_2}{\eta_1}\epsilon_2 + \dot{\epsilon}_2,$$
$$\dot{\epsilon} = \frac{2\mu_4}{\eta_3}\epsilon_4 + \dot{\epsilon}_4.$$

Eliminating $\epsilon_4$ from eqn. $(6.17)_a$ and eqn. $(6.19)_b$ leads to

(6.20) $$\dot{\epsilon} = \frac{\sigma}{\eta_3} + \frac{\dot{\sigma}}{2\mu_4} - \frac{2\mu_2}{\eta_3}\epsilon_2 - \frac{\mu_2}{\mu_4}\dot{\epsilon}_2.$$

Solving eqn. $(6.19)_a$ and eqn. (6.20) simultaneously, we get

(6.21)
$$\epsilon_2 = \frac{\left(1 + \frac{\mu_2}{\mu_4}\right)\dot{\epsilon} - \frac{\sigma}{\eta_3} - \frac{\dot{\sigma}}{2\mu_4}}{\frac{2\mu_2}{\eta_1}\left(\frac{\mu_2}{\mu_4} - \frac{\eta_1}{\eta_3}\right)},$$
$$\dot{\epsilon}_2 = \frac{\left(1 + \frac{\eta_1}{\eta_3}\right)\dot{\epsilon} - \frac{\sigma}{\eta_3} - \frac{\dot{\sigma}}{2\mu_4}}{\frac{\eta_1}{\eta_3} - \frac{\mu_2}{\mu_4}}.$$

Now, differentiating eqn. $(6.21)_a$ with respect to time and equating it to eqn. $(6.21)_b$, we get

(6.22) $$\sigma + \left(\frac{\eta_1}{2\mu_2} + \frac{\eta_3}{2\mu_4}\right)\dot{\sigma} + \frac{\eta_1\eta_3}{4\mu_3\mu_4}\ddot{\sigma} = (\eta_1 + \eta_3)\dot{\epsilon} + \frac{\eta_1\eta_3}{2\mu_2}\left(1 + \frac{\mu_2}{\mu_4}\right)\ddot{\epsilon}.$$

This is of the same form as eqn. (1.1) with

(6.23) $$p_1 = \frac{\eta_1}{2\mu_2} + \frac{\eta_3}{2\mu_4}, \quad p_2 = \frac{\eta_1\eta_3}{4\mu_3\mu_4}, \quad q_1 = \eta_1 + \eta_3, \quad q_2 = \frac{\eta_1\eta_3}{2\mu_2}\left(1 + \frac{\mu_2}{\mu_4}\right).$$

## 7. Final Remarks

We have shown four sets of energy storage and rate of dissipation which lead to four different three dimensional constitutive relations, which reduce in one dimension to the model developed by Burgers (eqn. (1.1)). Each of these three dimensional models can claim equal status as representing the three dimensional generalization of Burgers' model. We have chosen two natural configurations instead of one in all of these models. This is to incorporate two relaxation times possessed by Burgers-like fluid bodies. For example, in an asphalt concrete mixture (which has been shown to exhibit Burgers-like fluid behaviour), the aggregate matrix has a small relaxation time whereas the asphalt mortar matrix has relatively larger relaxation time (see [8]) and the choice of two natural configurations seems natural. It is possible that several other choices for the stored energy and the rate of dissipation could lead to the same one dimensional model due to Burgers. Interestingly, the structure of the three dimensional models that we have developed are quite distinct.

## 8. Appendix

Now, from eqn. (3.21) and eqn. (3.32) in [17]

$$(A.1) \qquad \frac{\partial f(\mathbf{L}_1)}{\partial \mathbf{D}_1} = \frac{1}{2}\left(\frac{\partial f}{\partial \mathbf{L}_1} + \left(\frac{\partial f}{\partial \mathbf{L}_1}\right)^T\right).$$

Hence, using eqn. (A.1)

$$(A.2) \qquad \begin{aligned} \frac{\partial}{\partial \mathbf{D}_1}\mathbf{B}_3.\left(\mathbf{F}_2\mathbf{L}_1\mathbf{F}_2^{-1} + \mathbf{F}_2^{-T}\mathbf{L}_1^T\mathbf{F}_2^T\right) &= \frac{\partial}{\partial \mathbf{D}_1}2\mathbf{B}_3.\mathbf{F}_2\mathbf{L}_1\mathbf{F}_2^{-1} \\ &= \frac{\partial}{\partial \mathbf{D}_1}2\mathbf{F}_2^T\mathbf{B}_3\mathbf{F}_2^{-T}.\mathbf{L}_1 \\ &= \mathbf{F}_2^T\mathbf{B}_3\mathbf{F}_2^{-T} + \mathbf{F}_2^{-1}\mathbf{B}_3\mathbf{F}_2. \end{aligned}$$

*Authors' address*: *Satish Karra, K. R. Rajagopal*, Department of Mechanical Engineering Texas A&M, University College Station, TX 77843-3123, USA, e-mail: `satkarra@tamu.edu`, `krajagopal@mengr.tamu.edu`.